\documentclass[12pt]{article}
\usepackage{amssymb,amsmath}
\begin{document}
\font\SY=msam10
\def\N{{\mathbb N}}
\def\Z{{\mathbb Z}}
\def\R{{\mathbb R}}
\def\C{{\mathbb C}}
\def\T{{\mathbb T}}
\def\L{{\cal L}}
\def\zp{\Z_+}
\def\dist{{\rm dist}\,}
\def\li{{\rm span}\,}
\def\square{\hbox{\SY\char"03}}
\def\epsilon{\varepsilon}
\def\phi{\varphi}
\def\kappa{\varkappa}
\def\theorem#1{\smallskip\noindent
{\scshape Theorem} {\bf #1}{\bf .}\hskip 8pt\sl}
\def\defin#1{\smallskip\noindent
{\scshape Definition} {\bf #1}{\bf .}\hskip 6pt}
\def\prop#1{\smallskip\noindent
{\scshape Proposition} {\bf #1}{\bf .}\hskip 8pt\sl}
\def\lemma#1{\smallskip\noindent
{\scshape Lemma} {\bf #1}{\bf .}\hskip 6pt\sl}
\def\cor#1{\smallskip\noindent
{\scshape Corollary} {\bf #1}{\bf .}\hskip 6pt\sl}
\def\quest#1{\smallskip\noindent
{\scshape Question} {\bf #1}{\bf .}\hskip 6pt\sl}
\def\epr{\smallskip\rm}
\def\rem#1{\smallskip\noindent
{\scshape Remark}  {\bf #1}.\hskip 6pt}
\def\wz{\thinspace}
\def\proof{P\wz r\wz o\wz o\wz f.\hskip 6pt}
\def\leq{\leqslant}
\def\geq{\geqslant}
\def\ssub#1#2{#1_{{}_{{\scriptstyle #2}}}}
\def\H{{\cal H}}
\def\B{{\cal B}}
\def\W{{\cal W}}
\def\slim{\mathop{\hbox{$\overline{\hbox{\rm lim}}$}}}
\def\ilim{\mathop{\hbox{$\underline{\hbox{\rm lim}}$}}}
\def\co{\text{\rm conv}}
\def\cov{\overline{\text{\rm conv}}}

\title{Corners of multidimensional numerical ranges}

\author{S.~Shkarin}

\date{}

\maketitle

\smallskip

\leftline{Queen's University Belfast, Department of Pure
Mathematics} \leftline{University road, Belfast BT7 1NN, UK}
\leftline{\bf e-mail: \tt s.shkarin@qub.ac.uk}

\rm\normalsize

\begin{abstract}
The $n$-dimensional numerical range of a densely defined linear
operator $T$ on a complex Hilbert space $\H$ is the set of vectors
in $\C^n$ of the form $(\langle Te_1,e_1\rangle,\dots,\langle
Te_n,e_n\rangle)$, where $e_1,\dots,e_n$ is an orthonormal system in
$\H$, consisting of vectors from the domain of $T$. We prove that
the components of every corner point of the $n$-dimensional
numerical range are eigenvalues of $T$.
\end{abstract}

\section{Introduction}

Throughout this paper $\H$ is a complex Hilbert space, $\C$ is the
field of complex numbers, $\R$ is the field of real numbers and $\N$
is the set of positive integers. A {\it densely defined linear
operator} on $\H$ is a linear operator taking values in $\H$ and
defined on a dense linear subspace $D_T$ called the {\it domain} of
$T$. We denote the set of densely defined linear operators on $\H$
by ${\cal L}(H)$. Symbol $L(\H)$ stands for the algebra of bounded
linear operators on $\H$.

\defin1 Let $T\in {\cal L}(H)$, $n\in\N$ and ${\cal E}(n,T)$ be the set
of orthonormal systems ${\bf e}=\{e_1,\dots,e_n\}$ in $\H$ such that
$e_j\in D_T$ for $1\leq j\leq n$. The {\it $n$-dimensional numerical
range} of $T$ is the set
$$
{\cal W}_n(T)=\{\tau(T,{\bf e}):{\bf e}\in {\cal E}(n,T)\},\ \
\text{where}\ \ \tau(T,{\bf e})=(\langle
Te_1,e_1\rangle,\dots,\langle Te_n,e_n\rangle)\in\C^n.
$$

Thus, ${\cal W}_n(T)$ is a subset of $\C^n$. Clearly the
$1$-dimensional numerical range of $T\in L(\H)$ is exactly the
conventional numerical range. This and related numerical range type
concepts were studied in
\cite{tsing1,tsing2,tsing3,hubner,mfil,bebi1,saf}. We would like to
mention the following elementary properties of ${\cal W}_n(T)$.

\begin{itemize}
\item[(E0)]${\cal W}_n(T)=\varnothing$ if
and only if ${\rm dim}\,\H<n$;
\item[(E1)]The set ${\cal W}_n(T)$ is permutation invariant. That
is, for any permutation $\pi\in S_n$, $U_\pi({\cal W}_n(T))={\cal
W}_n(T)$, where $U_\pi:\C^n\to\C^n$ is the linear map such that
$Ue_k=e_{\pi(k)}$, $\{e_1,\dots,e_n\}$ being the standard basis of
$\C^n$;
\item[(E2)]Let $T\in{\cal L}(\H)$, $\H_0$ be a closed linear subspace of $\H$
and $P$ be the orthoprojection onto $\H_0$. Suppose also that
$L=D_T\cap \H_0$ is dense in $\H_0$ and $T_0$ be the operator on
$\H_0$ with the domain $D_{T_0}=L$ defined by the formula
$T_0x=PTx$. Then ${\cal W}_n(T_0)\subseteq{\cal W}_n(T)$ for each
$n$;
\item[(E3)]Let $1\leq k\leq n$ and $\pi_k:\C^n\to\C^k$ be the
restriction operator: $\pi_k(z_1,\dots,z_n)=(z_1,\dots,z_k)$. Then
$\pi_k({\cal W}_n(T))={\cal W}_k(T)$ for any $T\in{\cal L}(\H)$ with
${\rm dim}\,\H\geq n$.
\end{itemize}

We always assume that $\C^n$ is equipped with the standard inner
product $\langle u,v\rangle=\sum\limits_{j=1}^n u_n\overline{v_n}$
and the $|\cdot|$ stands for the Euclidean norm: $|u|=\sqrt{\langle
u,u\rangle}$. We recall the definition of corner points.

\defin2 Let $\Lambda\subset\C^n$ and $u\in\Lambda$. Then $u$ is
called a {\it corner point of} $\Lambda$ or simply a {\it corner} of
$\Lambda$ if there exist $\epsilon>0$, $\delta>0$ and non-zero
$w\in\C^n$ such that
$$
\frac{\langle v-u,w\rangle}{|v-u|\cdot|w|}\geq\delta, \ \
\text{whenever $v\in\Lambda$ and $|v-u|<\epsilon$}.
$$

We are going to use the following slightly weaker property.

\defin3 Let $\Lambda\subset\C^n$ and $u\in\Lambda$. We say that $u$
is a {\it pseudocorner} of $\Lambda$ if there is no $\phi\in
C^1([-a,a],\C^n)$, $a>0$ such that $\phi([-a,a])\subseteq\Lambda$,
$\phi(0)=u$ and $\phi'(0)\neq 0$.

Clearly any pseudocorner is a corner. It is worth noting that $u$ is
not a pseudocorner of $\Lambda$ if and only if there exists a
diffeomorphism $F:\C^n\to\C^n$ such that $F(u)$ is the center of an
interval contained in $F(\Lambda)$.

\theorem{1.1}Let $T\in{\cal L}(\H)$, $n\in\N$,
$\lambda=(\lambda_1,\dots,\lambda_n)\in\C^n$ be a pseudocorner of
${\cal W}_n(T)$ and ${\bf e}\in {\cal E}(n,T)$ be such that
$\tau(T,{\bf e})=\lambda$. Then $Te_j=\lambda_j e_j$ for $1\leq
j\leq n$. In particular, $\lambda_j$ are eigenvalues of $T$. \epr

The Donoghue theorem \cite{don} says that any corner point of the
numerical range of $T\in L(\H)$ is an eigenvalue of $T$. This result
was strengthened by Sims \cite{sims}, who showed that the condition
to be a corner point of the numerical range can be replaced by the
weaker condition to be a boundary point of the numerical range such
that the curvature of the boundary in this point is infinite. One
can easily verify that such a point is a pseudocorner. This means
that Theorem~1.1 for the case $n=1$ and $T$ bounded is exactly the
Sims theorem. Thus, Theorem~1.1 generalizes the Donoghue and the
Sims theorems in two directions: for unbounded operators and for
multidimensional numerical ranges.

The case of corner points of the closure of the numerical range of
$T\in L(\H)$, which do not belong to the numerical range itself was
considered by H\"ubner \cite{hubner}, who proved that such points
are approximate eigenvalues of $T$. We generalize this result for
multidimensional numerical ranges.

\theorem{1.2}Let $T\in  L(\H)$, $n\in\N$,
$\lambda=(\lambda_1,\dots,\lambda_n)\in\C^n$ be a pseudocorner of
$\overline{{\cal W}_n(T)}$. Then $\lambda_j$ are approximate
eigenvalues of $T$. \epr

\section{Auxiliary lemmas}

\lemma{2.1}Let $X$ be a dense linear subspace of $\H$ and $Y$ be a
closed linear subspace of $\H$ such that ${\rm
dim}\,Y^\perp=n<\infty$. Then $X\cap Y$ is dense in $Y$. \epr

\proof Let $P$ be the orthoprojection in $\H$ onto $Y^\perp$. Since
$X$ is dense in $\H$, we see that $P(X)$ is dense in $Y^\perp$.
Taking into account that any dense linear subspace of a finite
dimensional Hausdorff topological vector space coincides with the
entire space, we have $P(X)=Y^\perp$. Therefore, we can choose an
$n$-dimensional linear subspace $L\subset X$ such that
$P(L)=Y^\perp$. Then $\H$ is the (non-orthogonal) direct sum of $L$
and Y. Since $L\subset X$, we see that $X$ is the direct sum of $L$
and $X\cap Y$. Since $L$ is finite dimensional, we obtain
$$
L\oplus Y=\H=\overline{X}=\overline{L\oplus X\cap Y}=L\oplus
\overline{X\cap Y}.
$$
This is possible only if $\overline{X\cap Y}=Y$. \square

\lemma{2.2}Let $T\in{\cal L}(\H)$, $n\in\N$,
$\lambda=(\lambda_1,\dots,\lambda_n)\in\C^n$ be a pseudocorner of
${\cal W}_n(T)$ and ${\bf e}\in {\cal E}(n,T)$ be such that
$\tau(T,{\bf e})=\lambda$. Let also $L=\li\{e_1,\dots,e_n\}$. Then
$T(L)\subseteq L$. \epr

\proof It suffices to show that $Te_j\in L$ for $1\leq j\leq n$. Let
$j\in\{1,\dots,n\}$. For any $u\in D_T\cap L^\perp$ with $\|u\|=1$
consider the map $\epsilon_j:[-1/2,1/2]\to \H$ defined by the
formula $\epsilon_j(t)=tu+(1-t^2)^{1/2}e_j$. For $t\in[-1/2,1/2]$
denote
$$
{\bf e}^t=\{e_1,\dots,e_{j-1},\epsilon_j(t),e_{j+1},\dots,e_n\}.
$$
Clearly ${\bf e}^t$ is an orthonormal system in $\H$, whose elements
belong to $D_T$. Thus, the map
$$
\phi:[-1/2,1/2]\to \C^n,\quad\phi(t)=\tau(T,{\bf e}^t)
$$
takes values in ${\cal W}_n(T)$. From the definition of $\phi$ it
follows that $\phi$ is infinitely differentiable. Indeed, any
$\phi_k$ for $k\neq j$ are constants and
$$
\phi_j(t)=\langle T\epsilon_j(t),\epsilon_j(t)\rangle= t^2 \langle
Tu,u\rangle+(1-t^2)\langle Te_j,e_j\rangle+t(1-t^2)^{1/2}(\langle
Tu,e_j\rangle+\langle Te_j,u\rangle).
$$
Differentiating the last display, we see
$$
\phi'_j(0)=\langle Tu,e_j\rangle+\langle Te_j,u\rangle.
$$
On the other hand since $\lambda=\phi(0)$ is a pseudocorner of
${\cal W}_n(T)$, we have $\phi'(0)=0$. Hence
\begin{equation}
\langle Tu,e_j\rangle+\langle Te_j,u\rangle=0\ \ \text{for each}\ \
u\in D_T\cap L^\perp. \label{f1}
\end{equation}
Substituting $iu$ instead of $u$ into (\ref{f1}), we obtain
\begin{equation}
\langle Tu,e_j\rangle-\langle Te_j,u\rangle=0\ \ \text{for each}\ \
u\in D_T\cap L^\perp. \label{f2}
\end{equation}
Subtracting (\ref{f2}) from (\ref{f1}), we get
\begin{equation}
\langle Te_j,u\rangle=0\ \ \text{for each}\ \ u\in D_T\cap L^\perp.
\label{f3}
\end{equation}
According to Lemma~2.1, $D_T\cap L^\perp$ is dense in $L^\perp$.
Therefore (\ref{f3}) implies that $\langle Te_j,u\rangle=0$ for each
$u\in L^\perp$. That is, $e_j\in (L^\perp)^\perp=L$. \square

The following lemma is an approximate version of Lemma~2.2.
Unfortunately it fails for general unbounded operators.

\lemma{2.3}Let $T\in L(\H)$, $n\in\N$,
$\lambda=(\lambda_1,\dots,\lambda_n)\in\C^n$ be a pseudocorner of
$\overline{{\cal W}_n(T)}$,
$\{\lambda^m=(\lambda^m_1,\dots,\lambda^m_n)\}_{m\in\N}$ be a
sequence of elements of ${\cal W}_n(T)$, converging to $\lambda$,
$\{{\bf e}^m\}_{m\in\N}$ be a sequence of elements of ${\cal
E}(n,T)$ such that $\tau(T,{\bf e^m})=\lambda^m$ and
$L_m=\li\{e_1^m,\dots,e_n^m\}$. Then ${\rm dist}\,(Te_j^m,L_m)\to 0$
as $m\to\infty$ for $1\leq j\leq n$. \epr

\proof Suppose the contrary. Then there exists $j\in\{1,\dots,n\}$
such that ${\rm dist}\,(Te_j^m,L_m)\not\to 0$. Passing to a
subsequence if necessary, we can assume that there exists
$\epsilon>0$ such that ${\rm dist}\,(Te_j^m,L_m)\geq\epsilon$ for
any $m\in\N$. Then we can choose $v_n\in L_m^\perp$ such that
$\|v_m\|=1$ and $|\langle Te_j^m,v_m\rangle|\geq \epsilon$ for each
$m\in\N$. Hence
\begin{equation}
\max\{|\langle Te_j^m,v_m\rangle+\langle Tv_m,e_j^m\rangle|,|\langle
Te_j^m,v_m\rangle-\langle Tv_m,e_j^m\rangle|\}\geq \epsilon\ \
\text{for each $m\in\N$.} \label{f4}
\end{equation}
Denote
$$
u_m=\left\{\begin{array}{ll}v_m&\text{if $|\langle
Te_j^m,v_m\rangle+\langle Tv_m,e_j^m\rangle|\geq\epsilon;$}\\
iv_m&\text{otherwise.}\end{array}\right.
$$
From (\ref{f4}) it follows that
\begin{equation}
|\langle Te_j^m,u_m\rangle+\langle Tu_m,e_j^m\rangle|\geq \epsilon\
\ \text{for each $m\in\N$.} \label{f5}
\end{equation}
For any $m\in\N$ consider the map $\epsilon^m_j:[-1/2,1/2]\to \H$
defined by the formula $\epsilon^m_j(t)=tu_m+(1-t^2)^{1/2}e_j$. For
$t\in[-1/2,1/2]$ denote
$$
{\bf
e}^{m,t}=\{e_1,\dots,e_{j-1},\epsilon^m_j(t),e_{j+1},\dots,e_n\}.
$$
Clearly ${\bf e}^{m,t}$ is an orthonormal system in $\H$. Thus, the
map
$$
\phi^m:[-1/2,1/2]\to \C^n,\quad\phi_m(t)=\tau(T,{\bf e}^{m,t})
$$
takes values in ${\cal W}_n(T)$. Taking into account that $T$ is
bounded, we see that for any $k=0,1,\dots$, there exists $c_k>0$
such that $|(\phi^m)^{(k)}(t)|\leq c_k$ for any $m\in\N$ and any
$t\in[-1/2,1/2]$. Since the Fr\'echet space ${\cal
E}=C^\infty([-1/2,1/2],\C^n)$ is a Fr\'echet--Montel space, see for
instance \cite{pit}, the sequence $\phi^m$ has a subsequence,
$\phi^{m_k}$, converging in the Fr\'echet space $\cal E$ to
$\phi\in\cal E$. Then $\phi$ takes values in $\overline{{\cal
W}_n(T)}$. As in the proof of the previous lemma,
$$
(\phi^m_j)'(0)=\langle Tu_m,e_j\rangle+\langle Te_j,u_m\rangle.
$$
According to (\ref{f5}), $|(\phi^m_j)'(0)|\geq\epsilon$ for each
$m\in\N$. Hence $|\phi_j'(0)|\geq\epsilon$. In particular
$\phi'(0)\neq0$. On the other hand $\phi^m(0)=\lambda^m\to\lambda$
and therefore $\phi(0)=\lambda$. This means that $\lambda$ is not a
pseudocorner of $\overline{{\cal W}_n(T)}$. We have arrived to a
contradiction. \square

The following lemma is Theorem~1.1 for the particular case when the
dimension of $\H$ equals $n$.

\lemma{2.4}Let ${\rm dim}\,\H=n\in\N$, $T\in L(\H)$,
$\lambda=(\lambda_1,\dots,\lambda_n)\in\C^n$ be a pseudocorner of
$\overline{{\cal W}_n(T)}$ and ${\bf e}\in {\cal E}(n,T)$ be such
that $\tau(T,{\bf e})=\lambda$. Then $Te_j=\lambda_j e_j$ for $1\leq
j\leq n$. \epr

\proof Since $\langle Te_j,e_j\rangle=\lambda_j$ for $1\leq j\leq
n$, it suffices to verify that $\langle Te_j,e_k\rangle=0$ for
$1\leq j,k\leq n$, $j\neq k$. Suppose the contrary. Then there exist
different $j,k\in\{1,\dots,n\}$ such that $\langle
Te_j,e_k\rangle=\neq 0$. Choose $\alpha\in\R$ such that
$$
|e^{-i\alpha}\langle Te_j,e_k\rangle+e^{i\alpha}\langle
Te_k,e_j\rangle|=|\langle Te_j,e_k\rangle| + |\langle
Te_k,e_j\rangle|.
$$
Then
\begin{equation}
|e^{-i\alpha}\langle Te_j,e_k\rangle+e^{i\alpha}\langle
Te_k,e_j\rangle|=c>0. \label{f6}
\end{equation}
For each $t\in[-1/2,1/2]$ consider the orthonormal basis ${\bf
e}^t=\{e^t_1,\dots,e^t_n\}$ in $\H$ defined by the formula
$$
e^t_l=\left\{\begin{array}{ll}e_l&\text{if $l\neq j$ and $l\neq
k$;}\\ (1-t^2)^{1/2}e_j+e^{i\alpha}te_k&\text{if $l=j$;}\\
(1-t^2)^{1/2}e_k-e^{-i\alpha}te_j&\text{if $l=k$.}\end{array}\right.
$$
Clearly, the function $\phi:[-1/2,1/2]\to\C^n$, $\phi(t)=\tau(T,{\bf
e}^t)$ is infinitely differentiable, takes values in ${\cal W}_n(T)$
and $\phi(0)=\lambda$. Differentiating with respect to $t$ the
expression $\langle Te_j^t,e_j^t\rangle$, we obtain
$$
\phi_j'(0)=e^{-i\alpha}\langle Te_j,e_k\rangle+e^{i\alpha}\langle
Te_k,e_j\rangle.
$$
According to (\ref{f6}), $\phi'_j(0)\neq 0$ and therefore
$\phi'(0)\neq 0$. Since $\phi(0)=\lambda$, the point $\lambda$ is
not a pseudocorner of ${\cal W}_n(T)$. This contradiction completes
the proof. \square

\section{Proof of Theorems 1.1 and 1.2}

Let $T\in{\cal L}(\H)$, $n\in\N$,
$\lambda=(\lambda_1,\dots,\lambda_n)\in\C^n$ be a pseudocorner of
${\cal W}_n(T)$, ${\bf e}\in {\cal E}(n,T)$ be such that
$\tau(T,{\bf e})=\lambda$ and $\H_0=\li\{e_1,\dots,e_n\}$. According
to Lemma~2.2, $T(\H_0)\subseteq \H_0$. Let $T_0\in L(\H_0)$ be the
restriction of $T$ to $\H_0$. Taking (E2) into account, we see that
${\cal W}_n(T_0)\subset {\cal W}_n(T)$. Since $\lambda=\tau(T,{\bf
e})=\tau(T_0,{\bf e})\in {\cal W}(T_0)$, we see that $\lambda$ is a
pseudocorner of ${\cal W}(T_0)$. From Lemma~2.4 it follows that
$Te_j=T_0e_j=\lambda_je_j$ for $1\leq j\leq n$. The proof of
Theorem~1.1 is complete.

Let $T\in  L(\H)$, $n\in\N$ and
$\lambda=(\lambda_1,\dots,\lambda_n)\in\C^n$ be a pseudocorner of
$\overline{{\cal W}_n(T)}$. Pick a sequence
$\lambda^m=(\lambda^m_1,\dots,\lambda^m_n)$ of elements of ${\cal
W}_n(T)$ converging to $\lambda$ and a sequence ${\bf
e}^m=(e^m_1,\dots,e^m_n)$ of elements of ${\cal E}(n,T)$ such that
$\tau(T,{\bf e}^m)=\lambda^m$ for each $m\in\N$. For any $m\in\N$,
let  $L_m=\li\{e^m_1,\dots,e^m_n\}$, $P_m$ be the orthoprojection
onto $L_m$ and $T_m\in L(L_m)$ be operators defined by the formulas
$T_mx=P_mTx$. For each $m\in\N$ let $\{t^m_{j,k}\}_{j,k=1}^n$ be the
matrix of $T_m$ with respect to the orthonormal basis ${\bf e}^m$.
Since $T$ is bounded there exists $c>0$ such that $\|T_m\|\leq c$
for any $m\in\N$. Therefore the set $\{t^m_{j,k}:m\in\N,\ 1\leq
j,k\leq n\}$ is bounded. Passing to a subsequence if necessary, we
can assume that $t^m_{j,k}\to t_{j,k}\in\C$ as $m\to\infty$ for
$1\leq j,k\leq n$. Consider the operator $T_0$ with the matrix
$\{t_{j,k}\}_{j,k=1}^n$ acting on $\H_0=\C^n$ and let ${\bf
e}=(e_1,\dots,e_n)$ be the standard basis of $\C^n$. Since
$t^m_{j,k}\to t_{j,k}\in\C$ as $m\to\infty$, we see that ${\cal
W}_n(T_0)\subseteq \overline{{\cal W}_n(T)}$ and $\tau(T_0,{\bf
e})=\lambda$. Since $\lambda$ is a pseudocorner of $\overline{{\cal
W}_n(T)}$, it is also a pseudocorner of ${\cal W}_n(T_0)$. According
to Lemma~2.4, the matrix $\{t_{j,k}\}_{j,k=1}^n$ is diagonal with
$t_{j,j}=\lambda_j$ for $1\leq j\leq n$. According to Lemma~2.3
${\rm dist}\,(e^m_j,L_m)\to 0$ as $m\to\infty$ for $1\leq j\leq n$.
Since additionally the sequence $\{t^m_{j,k}\}_{j,k=1}^n$ converges
to the diagonal matrix with the numbers $\lambda_j$ on the diagonal,
this means that
$\|Te^m_j-\lambda_je^m_j\|=\|T^me^m_j-\lambda_je^m_j\|\to 0$ as
$m\to\infty$ for $1\leq j\leq m$. It follows that $\lambda_j$ are
approximate eigenvalues of $T$. The proof of Theorem~1.2 is
complete.

\bigskip

{\bf Acknowledgements.} \ Supported by British Engineering and
Physical Research Council Grant GR/T25552/01.


\small
\begin{thebibliography}{99}

\itemsep=-2pt

\bibitem{don}W.~Donoghue, \it On the numerical range of a bounded
operator, \it Michigan Math. J. \bf 4\rm\ (1957), 261--263

\bibitem{tsing1}N.~Tsing, \it On the shape of generalized numerical
ranges, \rm Linear and Multilinear Algebra \bf10\rm\ (1981),
173--182

\bibitem{sims}B.~Sims, \it On a connection between the numerical
range and spectrum of an operator on a Hilbert space, \rm J. London
Math. Soc. \bf8\rm\ (1974), 57--59

\bibitem{tsing2}C.~Li and N.~Tsing, \it On the $k$th matrix
numerical range, \rm Linear and Multilinear Algebra \bf 28\rm\
(1991), 229--239

\bibitem{tsing3}W.~Cheung and N.~Tsing, \it The $C$-numerical range
of matrices is star-shaped, \rm Linear and Multilinear Algebra \bf
41\rm\ (1996), 245--250

\bibitem{hubner}M.~H\"ubner, \it Spectrum where the boundary is not
round, \rm Rocky Mountain J. Math. \bf25\rm\ (1995), 1351--1355

\bibitem{mfil}M.~Markus and I.~Filipenko, \it Nondifferentiable
boundary points of higher numerical range, \rm Linear Algebra Appl.
\bf21\rm\ (1978), 217--232

\bibitem{bebi1}N.~Bebiano, \it Nondifferentiable points of
$\partial W_c(A)$, \rm Linear and Multilinear Algebra \bf 19\rm\
(1986), 249--257

\bibitem{pit}H.~Sch\"afer, \it Topological Vector Spaces, \rm
Macmillan, New York, 1966

\bibitem{saf}Yu.~Safarov, \it Birkhoff's theorem and
multidimensional numerical range, \rm J. Funct. Anal. \bf222\rm\
(2005), 61--97

\end{thebibliography}
\end{document}